\documentclass{article}
\usepackage{amssymb,amsmath}
\usepackage{graphicx}

\def\qed{\hfill {\hbox{${\vcenter{\vbox{               %HOLLOW SQUARE
   \hrule height 0.4pt\hbox{\vrule width 0.4pt height 6pt
   \kern5pt\vrule width 0.4pt}\hrule height 0.4pt}}}$}}}

\def\tr{\triangleright}

\newtheorem{theorem}{Theorem}
\newtheorem{definition}{Definition}

\newtheorem{example}{Example}
\newtheorem{remark}[example]{Remark}

{\qed\par\medskip}

\date{}

\title{\Large \textbf{Birack shadow modules and their link invariants}}

\author{Sam Nelson\footnote{Email:\ knots@esotericka.org} \and 
Katie Pelland\footnote{Email:\ katie.pelland@gmail.com }} 

\begin{document}

\maketitle

\begin{abstract}
We introduce an associative algebra $\mathbb{Z}[X,S]$ associated to a 
birack shadow and define enhancements of the birack counting invariant
for classical knots and links via representations of $\mathbb{Z}[X,S]$ 
known as \textit{shadow modules}. We provide examples which demonstrate 
that the shadow module enhanced invariants are not determined by the 
Alexander polynomial or the unenhanced birack counting invariants.
\end{abstract}

\begin{center}
\parbox{5.5in}{\textsc{Keywords:} biracks, birack shadows, shadow algebra,
shadow modules, link invariants, enhancements of counting invariants

\smallskip

}\end{center}

\section{\large \textbf{Introduction}}\label{I}

In \cite{AG} an associative algebra was introduced arising from a finite
quandle $X$, known as the \textit{quandle algebra} $\mathbb{Z}[X]$, with
representations known as \textit{quandle modules}. In \cite{CEGS} and 
later \cite{HHNYZ}, quandle modules and rack modules were used to define
enhancements of the quandle and rack counting invariants.

In \cite{BN} the rack algebra was generalized to the case of finite biracks.
In this paper we generalize the birack algebra further to the case of 
\textit{birack shadows}, pairs $X,S$ where $X$ is a birack and $S$ is a set
with an action of $X$ satisfying certain properties. The \textit{shadow
algebra} of a birack shadow $X,S$ and its representations, known as
\textit{shadow modules}, are used to further enhance the birack shadow
counting invariant for classical knots and links.

The paper is organized as follows. In section \ref{R} we recall the basics
of biracks and birack shadows. In section \ref{S} we introduce the shadow 
algebra and shadow modules. In section \ref{E} we define the shadow module 
enhanced counting invariant and compute some examples, including examples 
which show that the invariant is strictly stronger than the unenhanced birack
shadow counting invariant and is not determined by the Alexander
polynomial. We conclude with a few questions for
future research in section \ref{Q}.

\section{\large \textbf{Biracks and Shadows}}\label{R}

We begin by recalling a definition (see \cite{FRS,FJK} or \cite{N} for more).

\begin{definition}\textup{
Let $X$ be a set. A \textit{birack} structure on $X$ is an invertible
map $B:X\times X\to X\times X$ satisfying the conditions
\begin{list}{}{}
\item[(i)] $B$ is \textit{sideways invertible}, that is, there exists a 
unique map $S:X\times X\to X\times X$ satisfying for all $x,y\in X$
\[S(B_1(x,y),x)=(B_2(x,y),y),\]
\item[(ii)] $B$ is \textit{diagonally invertible}, that is, the components
$(S^{\pm 1}\circ \Delta)_1$ and $(S^{\pm 1}\circ \Delta)_2$
of the compositions $S\circ \Delta$ and $S^{-1}\circ\Delta$
of the diagonal map $\Delta(x)=(x,x)$ with
sideways map and its inverse are bijections, and
\item[(iii)] $B$ is a solution to the \textit{set-theoretic Yang-Baxter 
equation:}
\[(B\times I)(I\times B)(B\times I)=
(I\times B)(B\times I)(I\times B).\]
\end{list}
We will occasionally find it convenient to abbreviate $B_1(x,y)=y^x$ and 
$B_2(x,y)=x_y$.
}\end{definition}

We interpret the birack operation $B(x,y)$ as taking labels $x,y\in X$
on the input semiarcs at a positive crossing and assigning labels to the 
output semiarcs. The inverse map $B^{-1}$ represents going through a
negative crossing, and the sideways map represents going through a positive
crossing sideways from left to right.

\[\includegraphics{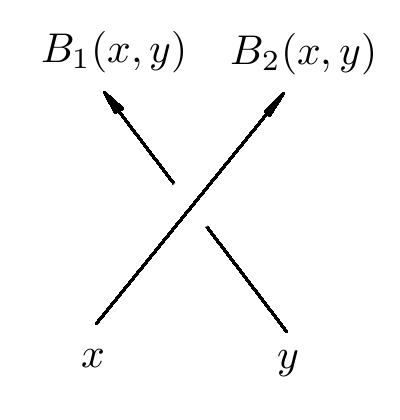} \quad
\includegraphics{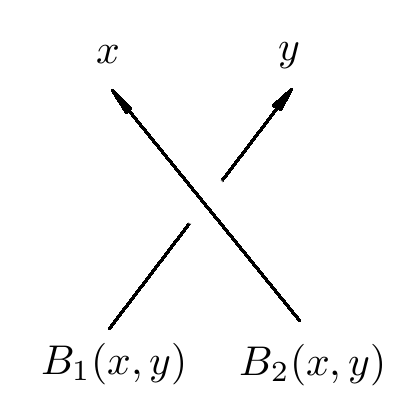} \quad
\includegraphics{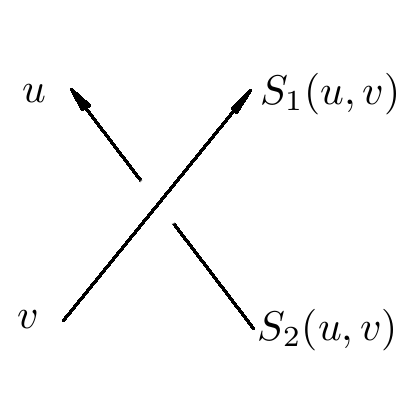}\]

The birack axioms are chosen such that labelings of the semiarcs in an
oriented blackboard-framed knot or link diagram are preserved by 
blackboard-framed isotopy moves. Invertibility and sideways invertibility
satisfy the direct and reverse type II moves, and the Yang-Baxter criterion
satisfies the type III move. Diagonal invertiblity satisfies the framed
type I moves.

\[\includegraphics{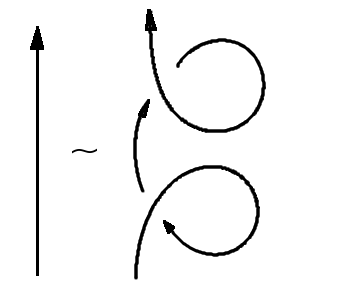} \quad
\includegraphics{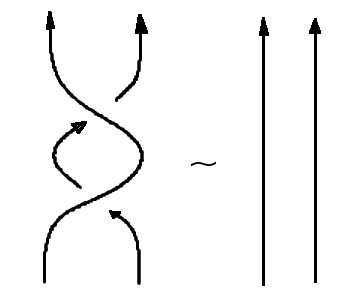} \quad
\includegraphics{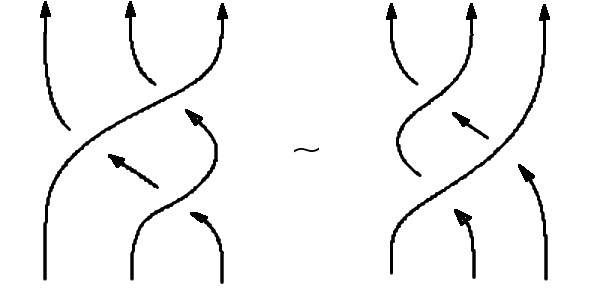} 
\]

The bijections $\alpha:X\to X$ and $\pi: X\to X$ defined by 
$\alpha=(S^{-1}\circ \Delta)_2^{-1}$ and $\pi=(S^{-1}\circ\Delta)_1\circ \alpha$
determine the labels on the semiarcs introduced in a blackboard framed
type I move.

\[\includegraphics{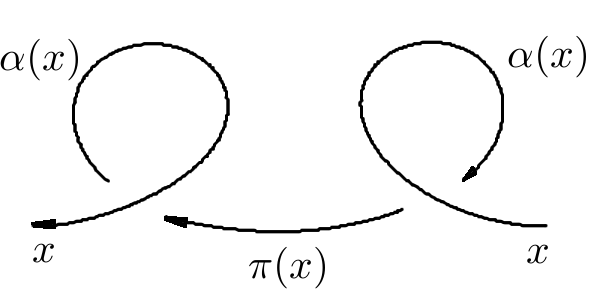}\]

If $X$ is a finite birack, then $\pi$ is an element of the finite symmetric
group $S_{|X|}$. The order or exponent of $\pi$, i.e. the smallest positive
integer $N$ such that $\pi^N(x)=x$ for all $x\in X$, is known as the 
\textit{birack rank} or \textit{birack characteristic} of $X$. If $X$ is 
a birack of rank $N$, then birack labelings of a link $L$ by $X$ are 
preserved by the $N$-phone cord move:

\[\includegraphics{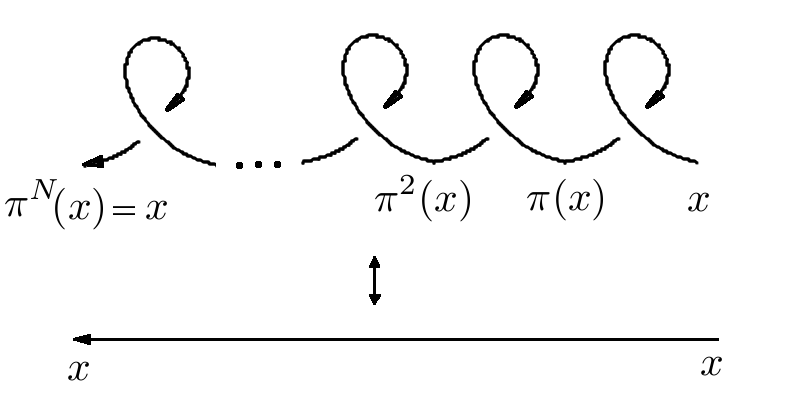}\]

Examples of biracks include 
\begin{itemize}
\item \textit{Groups.} A group $G$ is a birack under the map 
$B(x,y)=(x^nyx^{-n},x)$, for instance; many other birack structures on 
groups exist.
\item \textit{Quandles and racks.} A quandle $Q$ (see \cite{J,M}) or a rack
$R$ (see \cite{FR}) is a birack under
$B(x,y)=(y\tr x,x)$.
\item \textit{Biquandles.} A strong biquandle $X$ (see \cite{KR,FJK}) is a 
birack under $B(x,y)=(y^x,x_y)$.
\item \textit{Vector Spaces.} A vector space $V$ over a field $\mathbb{F}$
is a birack under $B(x,y)=(ty+sx,rx)$ where $t,s,r\in \mathbb{F}$ satisfy
$s^2=(1-tr)s$; biracks of this type are known as \textit{$(t,s,r)$-biracks}. 
See \cite{N}.
\end{itemize}

If $X=\{x_1,\dots,x_n\}$ is a finite birack, we can encode the birack structure
with a block matrix $M_X=[U|L]$ where the $i,j$ entry of $U$ is $k$ where
$x_k=B(x_j,x_i)$ and the $i,j$ entry of $L$ is $l$ where $x_l=B(x_i,x_j)$.
Note the reversed order of $i,j$ in $U$; this is for compatibility with
previous work.

We have the following standard notions:
\begin{definition}\textup{
Let $X$ and $X'$ be sets with birack structures $B$ and $B'$. Then we have:
\begin{itemize}
\item \textit{Homomorphisms.} A map $f:X\to X'$ is a 
\textit{birack homomorphism} if for all $x,y\in X$ we have
\[B'(f(x),f(y))=(f(B_1(x,y)),f(B_2(x,y))),\]
that is, if we have $B'\circ(f\times f)=(f\times f)\circ B$.
\item \textit{Subbiracks.} A subset $Y\subset X$ is a \textit{subbirack} of $X$ if the restriction
of $B$ to $Y\times Y$ defines a birack structure on $Y$; equivalently,
$Y\subset X$ is a subbirack if $Y\times Y$ is closed under $B$.
\end{itemize}
}\end{definition}

Given a finite birack $B$ of rank $N$ and a blackboard framed oriented link
diagram $L$ of $c$ components, the set of labelings of semiarcs in $L$ 
satisfying the crossing conditions
\[\includegraphics{sn-kp-6} \quad \includegraphics{sn-kp-7}\]
at every crossing corresponds bijectively with the set 
$\mathrm{Hom}(BR(L),B)$ of birack homomorphisms from the \textit{fundamental 
birack of $L$} (see \cite{FJK,N}) to $B$. Moreover, changing $L$ by blackboard 
framed moves and
$N$-phone cord moves also yields a bijection on the sets of labelings; 
thus, summing over a complete period of framings modulo $N$ yields a 
link invariant known as the \textit{integral birack counting invariant}
\[\Phi_{B}^{\mathbb{Z}}(L)=\sum_{\mathbf{w}\in(\mathbb{Z})^N} 
|\mathrm{Hom}(BR(L,\mathbf{w}),B)|\]
where $\mathbf{w}=(w_1,\dots, w_c)$ is the vector of framing numbers of $L$.
See \cite{N} for more.

If the elements of a birack are used to label the semiarcs of a an oriented
blackboard-framed classical link diagram on the sphere $S^2$, we can use 
elements of 
another set $S$ called \textit{shadows} to label the regions between the 
arcs. We can then let the birack label on a semiarc act on the label in one
region bounded by the semiarc to determine the label in the opposite
region as pictured below.

\[\includegraphics{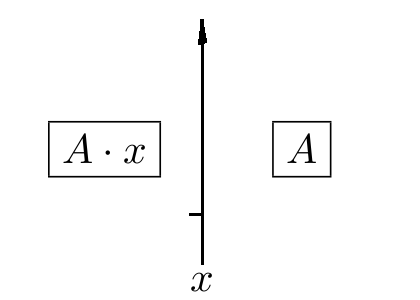}\]

We want birack shadows labelings to be preserved by blackboard-framed 
Reidemeister moves. To guarantee this, we need the shadow labelings
to be well-defined at crossings and at kinks.

\[\includegraphics{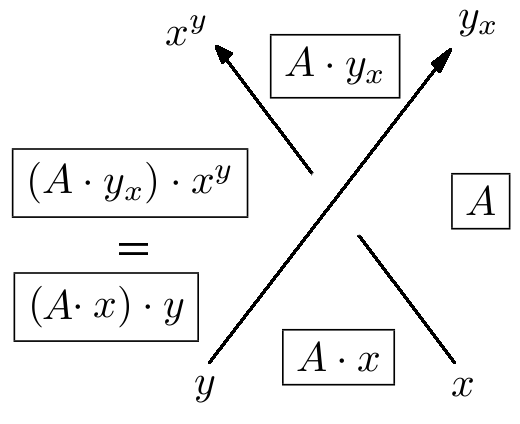} \quad \includegraphics{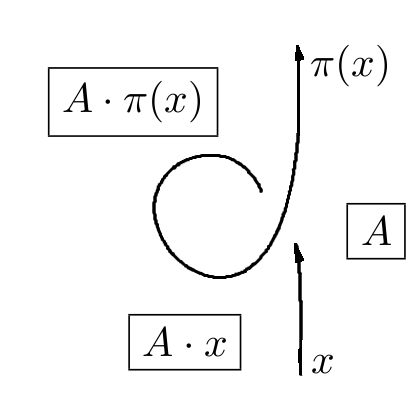}\]
Thus we have:
\begin{definition}
\textup{Let $X$ be a birack and $S$ a set. A \textit{birack shadow} structure
on $S$ is an invertible right action of $X$ on $S$ (i.e. a map 
$\cdot: X\times S\to S$)
satisfying for all $x,y\in X$ and $A\in S$,
\begin{list}{}{}
\item[(i)]{$(A\cdot y_x)\cdot x^y=(A\cdot x)\cdot y$ and}
\item[(ii)]{$A\cdot x=A\cdot\pi(x)$.}
\end{list}
We will refer to such an $S$ as an $X$-\textit{shadow}.}
\end{definition}

If $X=\{x_1,\dots,x_n\}$ is a finite birack and $S=\{A_1,\dots,A_m\}$ is a
finite set, we can encode an $X$-shadow structure on $S$ with an
$m\times n$ matrix $M_{X,S}$ whose $i,j$ entry is $k$ where $A_k=A_i\cdot x_j$.

Let $L$ be an oriented blackboard-framed link diagram on $S^2$, $X$ a finite
birack and $S$ a finite $X$-shadow. A \textit{shadow labeling} of $L$
is a labeling of the semiarcs of $L$ with elements of $X$ and the regions
between the semiarcs of $L$ with elements of $S$ such that at every crossing
and region boundary we have
\[\includegraphics{sn-kp-6} \quad 
\includegraphics{sn-kp-7} \quad 
\includegraphics{sn-kp-11} 
\]
While it is true that the number of shadow labelings 
$\Phi_{X,S}^{\mathbb{Z}}(L)$ of a link $L$ by a birack-shadow pair $(X,S)$ 
over a complete period of writhes is a link invariant, unfortunately 
$\Phi_{X,S}^{\mathbb{Z}}(L)$ is determined by the usual birack counting 
invariant $\Phi_{X}^{\mathbb{Z}}(L)$, since for any birack 
labeling we can choose a shadow label for a starting region and simply push
it across the semiarcs using the shadow operation and its inverse to determine
a unique valid shadow labeling. More precisely, we have

\begin{theorem}
The number of shadow labelings of a link $L$ by $(X,S)$ is given by
\[\Phi_{X,S}^{\mathbb{Z}}(L)=|S|\Phi_{X}^{\mathbb{Z}}(L).\]
\end{theorem}

Despite this fact, \textit{enhancements} of the shadow counting invariant
can give us information allowing us to distinguish links with equal birack
counting invariant. Previously studied examples include quandle 3-cocycle
invariants (where $X$ is a finite quandle, i.e. a rank $N=1$ birack with 
$B_2(x,y)=x$, $S=X$ and $A\cdot x=B_1(x,A)$, see \cite{CN2}) and rack shadow 
polynomials (see \cite{CN}). For the present paper, we will generalize the 
rack module idea from 
\cite{AG} to define \textit{shadow modules}, which we will use to enhance 
the birack shadow counting invariant.

\begin{remark}
\textup{If $L$ is a virtual link, the only shadow labelings which are 
preserved by virtual Reidemeister moves are constant labelings, i.e.
shadow labelings such that $A\cdot x=A$ for all $A\in S, x\in X$. 
Unfortunately, this fact cannot be used to detect non-classicality,
only non-planarity of a particular diagram of a virtual knot or link,
since even classical knots and links have non-planar virtual diagrams.} 
\end{remark}

\section{\large \textbf{The Shadow Algebra and Shadow Modules}}\label{S}

In \cite{AG} an associative algebra is defined with generators corresponding
to ordered pairs of elements of a finite quandle $X$, with relations obtained
from Reidemeister moves with a labeling of the strands by elements of $X$.
In \cite{HHNYZ} the quandle algebra was modified slightly for the purpose
of enhancing the integral rack counting invariant described in \cite{N1}.
In \cite{BN} the rack algebra was generalized to the cae of biracks.
We will now generalize one step further to the case of birack shadows.

Let $L$ be a classical knot or link diagram, $X$ a finite birack, $S$ an 
$X$-shadow and fix a shadow labeling of $L$ by $X,S$. We place 'beads' on the
semiarcs of $L$ related by linear equations with coefficients depending
on the rack and shadow labels at the crossings as depicted:

\[
\raisebox{-0.5in}{\includegraphics{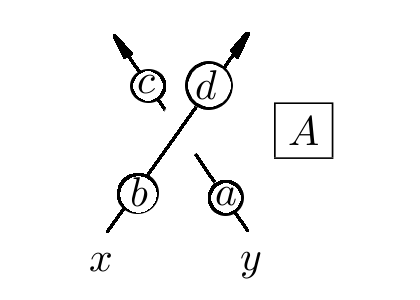}} \quad 
\raisebox{-0.5in}{\includegraphics{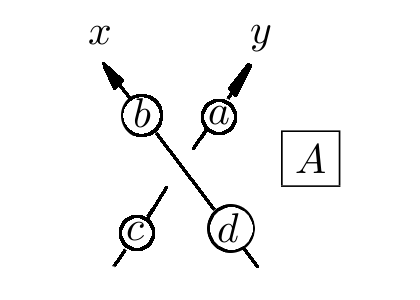}}\quad
\begin{array}{rcl}
c & = & t_{A,x,y} b+ s_{A,x,y}a \\
d & = & r_{A,x,y} a 
\end{array}
\]

\begin{definition}\label{sadef}
\textup{Let $X$ be a birack with birack rank $N$ and let 
$(S,\cdot)$ be a birack shadow. Let $\Omega[X,S]$ be the free 
$\mathbb{Z}$-algebra generated by elements of the form
$t_{A,x,y}^{\pm 1}, s_{A,x,y}$ and $r^{\pm 1}_{A,x,y}$
for $A\in S$, $x,y\in X$. Then the \textit{shadow algebra} 
$\mathbb{Z}[X,S]$ is the quotient of $\Omega[X,S]$ by the ideal $I$
generated by elements of the form}
\begin{itemize}
\item{$r_{A,x_{z^y},y_z}r_{A\cdot y_z,x,z^y}-              r_{A,x_y,z}r_{A\cdot z,x,y}$,}
\item{$t_{A,x_{z^y},y_z}r_{A,y,z}-         r_{A\cdot x_{yz},y^x,z^{x_y}}t_{A\cdot z,x,y}$,}
\item{$s_{A,x_{z^y},y_z}r_{A\cdot y_z,x,z^y}-r_{A\cdot x_{yz},y^x,z^{x_y}}s_{A\cdot z,x,y}$,}
\item{$t_{A\cdot y_z,x,z^y}t_{A,y,z}-t_{A\cdot x_{yz},y^x,z^{x_y}}t_{A,x_y,z}$,}
\item{$t_{A\cdot y_z,x,z^y}s_{A,y,z}-s_{A\cdot x_{yz},y^x,z^{x_y}}t_{A\cdot z,x,y}$,}
\item{$s_{A\cdot y_z,x,z^y}-        t_{A\cdot x_{yz},y^x,z^{x_y}}s_{A,x_y,z}r_{A\cdot z,x,y}
                                -s_{A\cdot x_{yz},y^x,z^{x_y}}s_{A\cdot z,x,y}$, 
\textup{and}}
\item{$1-\displaystyle{\prod_{k=0}^{N-1} 
(t_{A\cdot^{-1}\alpha(\pi^{k}(x)),\pi^k(x),\alpha(\pi^k(x))}
 r_{A\cdot^{-1}\alpha(\pi^{k}(x)),\pi^k(x),\alpha(\pi^k(x))}
+s_{A\cdot^{-1}\alpha(\pi^{k}(x)),\pi^k(x),\alpha(\pi^k(x))})}$.}
%\item{$t_{A,y_x,z_{x^y}}-t_{A\cdot x,y,z}$,}
%\item{$s_{A,y_x,z_{x^y}}-s_{A\cdot x,y,z}$,}
%\item{$t_{A\cdot z_{yx},x^{z_y},y^z}t_{A,x,z_y}-t_{A\cdot y_x,x^y,z}t_{A,x,y}$,}
%\item{$s_{A\cdot z_{yx},x^{z_y},y^z}t_{A\cdot x,y,z}-t_{A\cdot y_x,x^y,z}s_{A,x,y}$,}
%\item{$s_{A\cdot y_x,x^y,z}-t_{A\cdot z_{yx},x^{z_y},y^z}s_{A,x,z_y}-s_{A\cdot z_{yx},x^{z_y},y^z}s_{A\cdot x,y,z}$, \textup{and}}
%\item{$\displaystyle{1-\prod_{k=o}^{N-1} 
%(t_{A\cdot^{-1}\alpha(\pi^k(x)),\alpha(\pi^k(x)),\pi^{k}(x)}+
% s_{A\cdot^{-1}\alpha(\pi^k(x)),\alpha(\pi^k(x)),\pi^{k}(x)})}$.}
\end{itemize}
\textup{A \textit{shadow module} or $(X,S)$-module is a representation of 
$\mathbb{Z}[X,S]$, i.e. an abelian group $G$ with a family of automorphisms 
$t_{A,x,y},r_{A,x,y}:G\to G$
and endomorphisms $s_{A,x,y}:G\to G$ such that each of the above maps is
zero for all $A\in S$ and $x,y\in X$.}
\end{definition}

The conditions in definition \ref{sadef} come from comparing the beads
on the two sides
of the oriented framed Reidemeister moves and the $N$-phone cord move using 
the labeling convention defined above.
\[\includegraphics{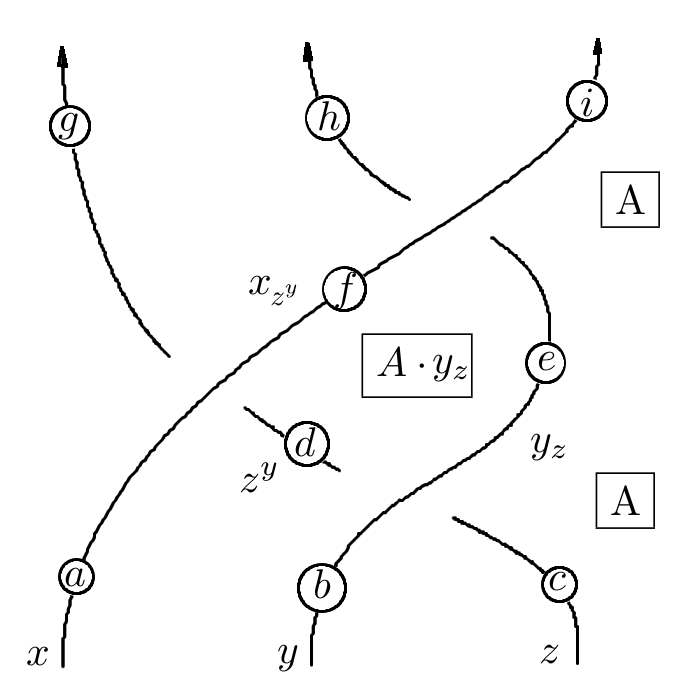}\quad \raisebox{1in}{$\sim$}\quad \includegraphics{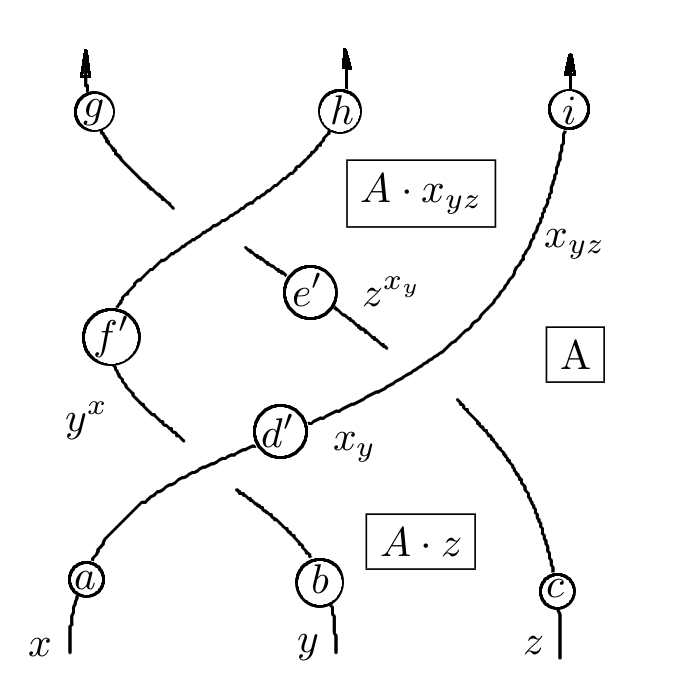}\]
For instance, comparing the coefficients of the bead $i$ before and after
the Reidemesiter III moves, we need 
$r_{A,x_{z^y},y_z}r_{A\cdot y_z,x,z^y}i=r_{A,x_y,z}r_{A\cdot z,x,y}i$.

Let $X=\{x_1,\dots, x_n\}$ be a finite birack, $S=\{A_1,\dots, A_m\}$ a 
finite $X$-shadow, and $k$ a ring with identity. We can define an 
$(X,S)$-module structure on $k$ by selecting $t_{A,x,y}, r_{A,x,y}\in k^{\times}$
and $s_{A,x,y}\in k$ such that the ideal in definition \ref{sadef} is zero;
the automorphisms and endomorphisms are then given by left multiplication by
the elements $t_{A,xy}, s_{A,x,y}$ and $r_{A,x,y}$. Such a module structure can
be conveniently encoded by an $3n\times kn$ block matrix
\[M_R=\left[\begin{array}{c|c|c}
T_1 & S_1 & R_1 \\ \hline
T_2 & S_2 & R_2 \\ \hline
\vdots & \vdots & \vdots \\ \hline
T_m & S_m & R_m \\ 
\end{array}\right]\]
where the $(i,j)$ entry of $T_l$ is $t_{A_l,x_i,x_j}$ and similarly for 
$S_l$ and $R_l$.

\begin{example}\label{exmod1}
\textup{Consider the birack $X$ and $X$-shadow $S$ with matrices}
\[M_X=\left[\begin{array}{cc|cc}
1 & 1 & 2 & 2 \\
2 & 2 & 1 & 1 \\
\end{array}\right],\quad
M_{X,S}=\left[\begin{array}{cc}
2 & 2 \\
3 & 3 \\
1 & 1
\end{array}\right].\]
\textup{Our \texttt{python} computations reveal 128 $X,S$-module structures
on $R=\mathbb{Z}_3$ including for instance}
\[M_{M}=\left[\begin{array}{cc|cc|cc}
1 & 1 & 2 & 2 & 2 & 2 \\
1 & 1 & 2 & 2 & 2 & 2 \\ \hline
2 & 1 & 2 & 1 & 1 & 1 \\
2 & 1 & 1 & 2 & 2 & 2 \\ \hline
2 & 1 & 2 & 2 & 2 & 2 \\
2 & 1 & 2 & 2 & 1 & 1 \\
\end{array}\right].\]
\end{example}

Let $X$ be a finite birack, $S$ an $X$-shadow, $L$ an oriented link
of $c$ components, and let $f$ be an $(X,S)$-labeling of a diagram of $L$.
Placing beads $a_1,\dots, a_n$ on each semiarc in $L$, we get a free 
$\mathbb{Z}[X,S]$-module generated by the beads. Each crossing gives
us two equations; the quotient of this free $\mathbb{Z}[X,S]$-module
by the crossing bead relations is then by construction aninvariant of
the $(X,S)$-labeling $f$ of $L$ under $(X,S)$-labled framed Reidemeister
moves and $N$-phone cord moves, called the \textit{fundamental} 
$\mathbb{Z}[X,S]$-\textit{module of} $f$, denoted $\mathbb{Z}_{X,S}[f]$.
We can express $\mathbb{Z}_{X,S}[f]$ with a coefficient matrix of the
homogeneous system given by the crossing equations.

\begin{example}
\textup{Let $X,S$ be as in example \ref{exmod1} and consider the labeling
of the figure eight knot $4_1$ pictured below. We obtain the listed 
presentation matrix for $\mathbb{Z}_{X,S}[f]$.}
\[\raisebox{-0.5in}{\includegraphics{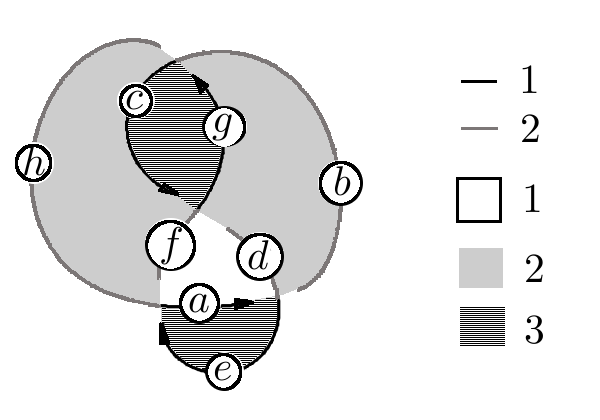}}
\scalebox{0.9}{$\quad M_{\mathbb{Z}_{X,S}[f]}=
\left[
\begin{array}{cccccccc}
0 & 0 & 0 & 0 & t_{3,2,1} & -1 & 0 & s_{3,2,1} \\
-1 & 0 & 0 & 0 & 0 & 0 & 0 & r_{3,2,1} \\
t_{3,2,1}& -1 & 0 & s_{3,2,1} & 0 & 0 & 0 & 0 \\
0 & 0 & 0 & r_{3,2,1} & -1 & 0 & 0 & 0 \\
0 & 0 & s_{1,1,2} & 0 & 0 & 0 & -1 & t_{1,1,2} \\
0 & -1 & r_{1,1,2} & 0 & 0 & 0 & 0 & 0 \\
0 & 0 & -1 & t_{1,1,2} & 0 & 0 & s_{1,1,2} & 0 \\
0 & 0 & 0 & 0 & 0 & -1 & r_{1,1,2} & 0 \\   
\end{array}
\right]$}
\] 
\end{example}

\section{\large \textbf{Enhancing the Counting Invariant}}\label{E}

Let $L$ be a blackboard framed oriented link, $X$ a birack of finite rank
$N$, $S$ a finite $X$-shadow and $M$ a $\mathbb{Z}[X,S]$-module. 

\begin{definition}\textup{
The \textit{shadow module multiset} invariant of a link $L$ of $c$ components
with respect to the $X,S$-shadow module $M$ is the multiset
\[\Phi^{M}_{X,S,M}(L)=\{|\mathrm{Hom}_{\mathbb{Z}[X,S]}(\mathbb{Z}_{X,S}[f],M)|\ :
\ f\in \mathcal{L}((L,\mathbf{w}),(X,S)),\ \mathbf{w}\in(\mathbb{Z}_N)^c
\}\]
where $\mathcal{L}((L,\mathbf{w}),(X,S))$ is the set of shadow labelings
of the link diagram $L$ with writhe vector $\mathbf{w}=(w_1,\dots,w_c)$.}

\textup{The \textit{shadow module polynomial} invariant is
\[\Phi_{X,S,M}(L)=\sum_{\mathbf{w}\in(\mathbb{Z}_N)^c} 
\left(\sum_{f\in \mathcal{L}((L,\mathbf{w}),(X,S))}
u^{|\mathrm{Hom}_{\mathbb{Z}[X,S]}(\mathbb{Z}_{X,S}[f],M)|}\right).\]
}\end{definition}

By construction, for every bead labeling of an $X,S$-labeling of $L$ by
beads in $M$, there is a unique corresponding bead labeling of every
$X,S$-labeled diagram obtained from $L$ by blackboard framed oriented
Reidemeister moves and $N$-phone cord moves. Hence, the set of such 
$M$-labelings is a signature of the $X,S$-labeling of $L$, and the 
multiset of these signatures over the set of all $X,S$-labelings forms 
an enhancement of the birack counting invariant. More formally, we have:

\begin{theorem}
If $L$ and $L'$ are ambient isotopic classical links, then for any
finite birack $X$, $X$-shadow $S$ and $X,S$-module $M$, we have
\[\Phi^{M}_{X,S,M}(L) =\Phi^{M}_{X,S,M}(L')\quad \mathit{and}\quad
\Phi_{X,S,M}(L) =\Phi_{X,S,M}(L').\]
\end{theorem}

\begin{remark}\textup{The birack module invariants defined in \cite{BN}
coincide with the special case of $\Phi_{X,S}^{M}$ where the shadow $S$ is a
singleton set for classical links.}
\end{remark}

To compute the shadow module invariant of a link $L$ with respect to 
an $X,S$-module structure $M$ on a commutative ring $G$, we obtain a 
presentation matrix for $\mathbb{Z}_{X,S}[f]$ for each shadow labeling $f$ 
of $L$ and replace the $t_{A,x,y}, s_{A,x,y}$ and $r_{A,x,y}$ with their 
values in $M$. The cardinality $n$ of the solution space of resulting matrix 
is then the signature of $f$, and $f$ contributes $u^n$ to the overall
invariant. We then repeat this computation for all shadow labelings over all
writhe vectors of $L$.

\begin{example}
\textup{For a simple example, let us compute $\Phi_{X,S}^{M}(K)$ for the 
trefoil knot $3_1$ with respect to the birack, shadow, and $X,S$-module 
structure on  $\mathbb{Z}_3$ from example \ref{exmod1}.
Since the birack rank of $X$ is $N=2$, we need to consider diagrams of
$3_1$ with writhes equal to $0$ and $1$ mod $2$. There are six total 
$X,S$-labelings of an even-writhe diagram of $3_1$ and no valid labelings 
of an odd-writhe diagram. The labeling $f$ below has the listed 
presentation matrix for $\mathbb{Z}_{X,S}[f]$.}
\[\includegraphics{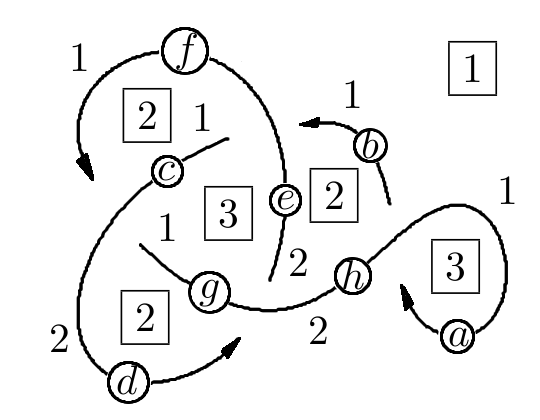}\quad\raisebox{0.55in}{\scalebox{0.9}{
$M_{\mathbb{Z}_{X,S}[f]}
=\left[\begin{array}{cccccccc}
t_{3,2,1} & -1 & 0 & 0 & 0 & 0 & 0 & s_{3,2,1} \\
-1 & 0 & 0 & 0 & 0 & 0 & 0 & r_{3,2,1} \\
0 & t_{1,2,1} & -1 & 0 & s_{1,2,1} & 0 & 0 & 0 \\
0 & 0 & 0 & 0 & r_{1,2,1} & -1 & 0 & 0 \\
0 & 0 & s_{1,1,1} & 0 & 0 & t_{1,1,1} & -1 & 0 \\
0 & 0 & r_{1,1,1} & -1 & 0 & 0 & 0 & 0 \\
0 & 0 & 0 & t_{1,1,2} & -1 & 0 & s_{1,1,2} & 0 \\
0 & 0 & 0 & 0 & 0 & 0 & r_{1,1,2} & -1 \\
\end{array}\right]$}}\]
\textup{Replacing $t_{A,x,y}, s_{A,x,y}$ and $r_{A,x,y}$ with their 
values in $R$ and row-reducing over $\mathbb{Z}_3$, we have}
\[
\left[\begin{array}{cccccccc}
2 & 2 & 0 & 0 & 0 & 0 & 0 & 2 \\
2 & 0 & 0 & 0 & 0 & 0 & 0 & 1 \\
0 & 1 & 2 & 0 & 2 & 0 & 0 & 0 \\
0 & 0 & 0 & 0 & 2 & 2 & 0 & 0 \\
0 & 0 & 2 & 0 & 0 & 1 & 2 & 0 \\
0 & 0 & 2 & 2 & 0 & 0 & 0 & 0 \\
0 & 0 & 0 & 1 & 2 & 0 & 2 & 0 \\
0 & 0 & 0 & 0 & 0 & 0 & 2 & 2 \\
\end{array}\right] \rightarrow
\left[\begin{array}{cccccccc}
1 & 1 & 0 & 0 & 0 & 0 & 0 & 1 \\
0 & 1 & 0 & 0 & 0 & 0 & 0 & 2 \\
0 & 0 & 1 & 0 & 1 & 0 & 0 & 2 \\
0 & 0 & 0 & 1 & 2 & 0 & 0 & 1 \\
0 & 0 & 0 & 0 & 1 & 2 & 2 & 2 \\
0 & 0 & 0 & 0 & 0 & 0 & 1 & 1 \\
0 & 0 & 0 & 0 & 0 & 0 & 0 & 0 \\
0 & 0 & 0 & 0 & 0 & 0 & 0 & 0 \\
\end{array}\right]
\]
\textup{Thus, the shadow labeling $f$ contributes $u^9$ to 
$\Phi_{X,S}^{M}(3_1)$. Repeating this for all shadow labelings, we obtain
$\Phi_{X,S}^{M}(3_1)=6u^9$; this detects the knottedness of $3_1$, since
the unknot has $\Phi_{X,S}^{M}(\mathrm{Unknot})=6u^3$. We note that this
shows that $\Phi_{X,S}^{M}$ is a properly stronger invariant than the 
unenhanced shadow counting invariant, since all classical knots have
shadow counting invariant value $\Phi_{X,S}^{\mathbb{Z}}=6$.}
\end{example}

\begin{example}
\textup{Let $X,S$ and $M$ be the birack, shadow, and $X,S$-module structure on 
$\mathbb{Z}_5$ with matrices
\[M_B=\left[\begin{array}{cc|cc}
1 & 1 & 2 & 2 \\
2 & 2 & 1 & 1 \\
\end{array}\right],\quad 
M_s=\left[\begin{array}{cc}
2 & 2 \\
1 & 1
\end{array}\right],\quad 
M_{R}=\left[\begin{array}{cc|cc|cc}
1 & 2 & 2 & 2 & 2 & 2 \\ 
1 & 2 & 2 & 2 & 4 & 4 \\ \hline
1 & 3 & 2 & 1 & 4 & 4 \\
1 & 3 & 4 & 2 & 3 & 3 \\
\end{array}\right]
\]
We computed $\Phi_{X,S}^{R}$ for all prime classical
knots with up to eight crossings and all prime classical links with up to
seven crossings; these are listed in the table below.}
\[\begin{array}{l|l}
\Phi_{X,S}^{R} & L \\ \hline
4u^5 & 3_1, 5_2,6_1,6_2,6_3,7_1,7_2,7_3,7_5,7_6,7_7,8_1,8_2,8_3,8_4,8_5,8_6,8_7,8_{10}, 8_{12}, 8_{13}, 8_{14}, 8_{15},8_{17},8_{19}, 8_{20} \\
4u^{25} & 4_1, 5_1, 7_4, 8_8,8_9, 8_{11}, 8_{16}, 8_{18}, 8_{21} \\ 
8u^5 & L2a1, L4a1, L5a1, L6a1,L6a3, L7a1,L7a3,L7a4,L7a5,L7a6,L7n1, L7n2 \\
8u^{25} & L6a2,L7a2 \\
16u^5 & L6a4, L6a5,L6n1 \\
16u^{25} & L7a7\\
\end{array}
\]
\end{example}

\begin{example}
\textup{For our final example, we note that our \texttt{python} computations
show that $\Phi_{X,S}^{R}$ is not determined by the Alexander
polynomial. The $X,S$-module structure $M$ on $\mathbb{Z}_3$ given by
\[M_X=\left[\begin{array}{ccc|ccc}
1 & 3 & 1 & 3 & 3 & 3 \\
2 & 2 & 2 & 2 & 2 & 2 \\
3 & 1 & 3 & 1 & 1 & 1 \\
\end{array}\right],\ 
M_S=\left[\begin{array}{ccc}
2 & 2 & 2 \\
1 & 1 & 1
\end{array}\right],\
M_{M}=\left[\begin{array}{ccc|ccc|ccc}
1 & 1 & 1 & 2 & 0 & 1 & 1 & 1 & 1 \\
1 & 1 & 1 & 0 & 0 & 0 & 2 & 1 & 2 \\
1 & 1 & 1 & 1 & 0 & 2 & 1 & 1 & 1 \\ \hline
2 & 2 & 2 & 2 & 0 & 1 & 2 & 2 & 2 \\
2 & 2 & 2 & 0 & 0 & 0 & 1 & 2 & 1 \\
2 & 2 & 2 & 1 & 0 & 2 & 2 & 2 & 2 \\
\end{array}\right]
\]
detects the difference between the Alexander-equivalent knots $8_{18}$ and 
$9_{24}$ with $\Phi_{X,S}^{R}(8_{18})=4u^3+4u^{27}\ne 4u^3+4u^9
=\Phi_{X,S}^{R}(9_{24})$. Our \texttt{python} code computed the values
of $\Phi_{X,S}^{R}$ for the prime knots with up to eight crossings
in the table below.}
\[\begin{array}{l|l}
\Phi_{X,S}^{R} & L \\ \hline
4u^3+4u^9 & 3_1, 6_1, 7_4, 7_7,8_5,8_{10}, 8_{11}, 8_{15}, 8_{19}, 8_{20},8_{21} \\ 
4u^3+4u^{27} & 8_{18} \\
8u^3 & 4_1, 5_1, 5_2, 6_2, 6_3, 7_1, 7_2, 7_3, 7_5, 7_6,8_1,8_2,8_3,8_4,8_6,
8_7,8_8,8_9,8_{12}, 8_{13}, 8_{14}, 8_{16}, 8_{17} \\
\end{array}\]
\end{example}

\section{\large \textbf{Questions for Future Research}}\label{Q}

In this section we collect a few questions for future research.

Quandle and biquandle labelings have been used with 3-cocycles to define
invariants of knotted links and surfaces in $R^4$.
What is the correct generalization of shadow modules to the case of
knotted surfaces?

We have defined only the simplest possible invariant using shadow modules,
namely counting shadow module labelings. Any quantity computable from a
shadow module labeled diagram that is preserved by labeled Reidemeister moves
defines and enhancement of the counting invariant.  Obvious ideas include
rack/birack homology with shadow module coefficients and homology theories
with chain groups generated by shadow module elements themselves. What 
enhancements of the shadow module counting invariant can be found?

Moreover, we have only computed examples of shadow module structures on
$\mathbb{Z}_n$; shadow modules on finite non-commutative rings such as
$(\mathbb{Z}_n)^m$ should prove interesting.

We have found examples showing that $\Phi_{X,S}^{M}$ is not determined by the
Alexander polynomial. What about other invariants such as twisted Alexander 
polynomials, Jones and colored Jones polynomials, HOMFLYpt, Khovanov and 
Knot Floer homologies?

\medskip

\texttt{python} code for computing the invariants described in this paper
is available from the second author's website at 
\texttt{http://www.esotericka.org}.

\bigskip

\noindent\textsc{Department of Mathematics \\ 
Claremont McKenna College \\
850 Columbia Ave. \\
Claremont, CA 91711} \\
\texttt{knots@esotericka.org}

\medskip

\noindent\textsc{Department of Mathematics \\ 
Pomona College \\
610 North College Ave \\
Claremont, CA 91711} \\
\texttt{}


\begin{thebibliography}{10}

\bibitem{AG}{N. Andruskiewitsch and M. Gra\~{n}a.
From racks to pointed Hopf algebras.
\textit{Adv. Math.} \textbf{178} (2003) 177-243.}

\bibitem{KA}{D. Bar-Natan (Ed.). The Knot Atlas.
\texttt{http://katlas.math.toronto.edu/wiki/Main\underline{\ }Page}}

\bibitem{BN}{R. Bauernschmidt and S. Nelson. Birack modules and their
link invariants. arXiv:1103.0301 }

\bibitem{CN}{T. Carrell and S. Nelson. On rack polynomials. arXiv:0809.5075,
to appear in \textit{J. Alg. Appl.}}

\bibitem{CEGS}{J. S. Carter, M. Elhamdadi, M. Gra\~na and
M. Saito. Cocycle knot invariants from quandle modules and 
generalized quandle homology.
\textit{Osaka J. Math.} \textbf{42} (2005) 499-541.}

\bibitem{CN2}{W. Chang and S. Nelson. Rack shadows and their invariants.  
arXiv:0910.3002. To appear in  \textit{J.\ Knot Theory Ramifications.}}

\bibitem{FJK}{R. Fenn, M. Jordan-Santana and L. Kauffman. Biquandles 
and virtual links.  \textit{Topology Appl.}  \textbf{145}  (2004) 157-175.}

\bibitem{FR}{R. Fenn and C. Rourke.
 Racks and links in codimension two.
 \textit{J. Knot Theory Ramifications}  \textbf{1}  (1992) 343-406.}

\bibitem{FRS}{R. Fenn, C. Rourke and B. Sanderson. 
Trunks and classifying spaces. \textit{Appl. Categ. Structures} \textbf{3} 
(1995) 321--356.}

\bibitem{HHNYZ}{A. Haas, G. Heckel, S. Nelson, J. Yuen, Q. Zhang. 
Rack Module Enhancements of Counting Invariants. arXiv:1008.0114, to appear in
\textit{Osaka J. Math.}}

\bibitem{J}{D. Joyce.
 A classifying invariant of knots, the knot quandle.
 \textit{J. Pure Appl. Algebra}  \textbf{23}  (1982)  37-65.}


\bibitem{KR}{L. H. Kauffman and D. Radford. Bi-oriented quantum algebras, 
and a generalized Alexander polynomial for virtual links. 
\textit{Contemp. Math}. \textbf{318} (2003) 113-140.}

\bibitem{M}{S. V. Matveev.
Distributive groupoids in knot theory.
\textit{Math. USSR, Sb.} \textbf{47} (1984) 73-83.}

\bibitem{N1}{S. Nelson. Link invariants from finite racks. arXiv:0808.0029;
to appear in \textit{Fund. Math.} }

\bibitem{N}{S. Nelson. Link invariants from finite biracks.
arXiv:1002.3842.}


\end{thebibliography}
\end{document}